\documentclass[10pt,twocolumn,twoside,journal]{IEEEtran}

\usepackage{cite}
\usepackage{fixltx2e}
\usepackage{graphicx}
\usepackage[cmex10]{amsmath}
\usepackage{amsfonts} 
\usepackage{amssymb}

\def\R{{\mathbb{R}}}

\newtheorem{proposition}{Proposition}

\begin{document}

\title{Empirical wavelet transform}
\author{\thanks{Manuscript received October, 2012. Revised version received February, 2013. Copyright (c) 2012 IEEE. Personal use of this material is 
permitted. However, permission to use this material for any other purposes must be obtained from the IEEE by sending a request to pubs-permissions@ieee.org} 
J\'er\^ome Gilles\thanks{J. Gilles is with the Department of Mathematics, University of California, Los Angeles (UCLA), 520 Portola Plaza, 
Los Angeles, CA 90024, USA email: jegilles@math.ucla.edu.}
\thanks{This work was partially founded by the following grants NSF DMS-0914856, NSF DMS-1118971, ONR N00014-08-1-1119, ONR N0014-09-1-0360, ONR MURI USC, the UC Lab Fees Research and the Keck Foundation.}}

\markboth{IEEE Trans. on Signal Processing,~Vol.~xx, No.~xx, February~2013}%
{Gilles: Empirical wavelet transform}

\maketitle

\begin{abstract}
Some recent methods, like the Empirical Mode Decomposition (EMD), propose to decompose a signal accordingly to its contained information. Even though 
its adaptability seems useful for many applications, the main issue with this approach is its lack of theory.
This paper presents a new approach to build adaptive wavelets. The main idea is to extract the different modes of a signal by designing an appropriate 
wavelet filter bank. This construction leads us to a new wavelet transform, called the empirical wavelet transform. Many experiments are presented showing 
the usefulness of this method compared to the classic EMD.

\end{abstract}


\IEEEpeerreviewmaketitle

\section{Introduction}\label{sec:intro}
\IEEEPARstart{A}{daptive} methods to analyze a signal is of great interest to find sparse representations in the context of compressive sensing. ``Rigid'' methods, like the Fourier or wavelets transforms, correspond 
to the use of some basis (or frame) designed independently of the processed signal. The aim of adaptive methods is to construct such a basis directly based on the information contained in the signal. 
A well known way to build an adaptive representations is the basis pursuit approach which is used in the wavelet packets transform. Even though the wavelet packets have shown interesting results for 
practical applications, they still are based on a prescribed subdivision scheme. A completely different approach to build an adaptive representation is the algorithm called ``Empirical Mode Decomposition'' (EMD) proposed by Huang et al. 
\cite{Huang1998}. The purpose of this method is to detect the principal ``modes'' which represent the signal (roughly speaking, a mode corresponds to a signal which have a compactly supported Fourier spectrum). 
This method has gained a lot of interest in signal analysis this last decade, mainly because it is able to separate stationary and non-stationary components from a signal. 
However, the main issue of the EMD approach is its lack of mathematical theory. Indeed, it is an algorithmic approach and, due to its non-linearity, is 
difficult to model. Nevertheless, some experiments \cite{Flandrin2004a,Flandrin2005,Flandrin2004} show that EMD behaves like an adaptive filter bank.\\
Some recent works attempt to model EMD in a 
variational framework. In \cite{Daubechies2011}, the authors proposed to model a mode as an amplitude modulated-frequency modulated (AM-FM) signal and then use the properties of such signals to build a functional to
represent the whole signal. Then they are able to retrieve the different modes by minimizing this functional. Another proposed variational approach is the work of Hou et al. \cite{Hou2011} where the authors 
also use the AM-FM formalism. They propose to minimize a functional which is build on some regularity assumptions about the different components and uses higher-order total variation priors.\\
In this paper, we propose a new approach to build adaptive wavelets capable of extracting AM-FM components of a signal. The key idea is that such AM-FM components have a compact support Fourier spectrum. 
Separating the different modes is equivalent to segment the Fourier spectrum and to apply some filtering corresponding to each detected support. We will show that it is possible to adapt the wavelet 
formalism by considering distinct Fourier supports and then build a set of functions which form an orthonormal basis. Based on this construction, we propose an empirical wavelet transform (and its inverse) to analyze 
a signal.\\
The remainder of the paper is organized as follows. Section~\ref{sec:exist} has two distinct subsections: in \ref{sec:emd}, we recall the principle of the EMD algorithm and the AM-FM model; while in \ref{sec:wav}, 
we recall some wavelet formalism which will be useful in our own construction and we discuss some of the existing adaptive wavelet methods. In section~\ref{sec:ewt}, we build the proposed empirical wavelets and give some 
of their properties, then the empirical wavelet transform and its inverse are introduced. Section~\ref{sec:expe} show many experiments based on simulated and real signals. 
The time-frequency representation based on the Hilbert transform is introduced in section~\ref{sec:tfr}.
In section~\ref{sec:detect}, we address the question of the estimation of the number of modes. An extension to 2D signals (images) is presented in section~\ref{sec:2d}. 
Finally, we conclude and give some perspectives in section~\ref{sec:conc}.

\section{Existing approaches}\label{sec:exist}

\subsection{Empirical Mode Decomposition}\label{sec:emd}
\begin{figure}[!t]
\centering
\includegraphics[width=\columnwidth]{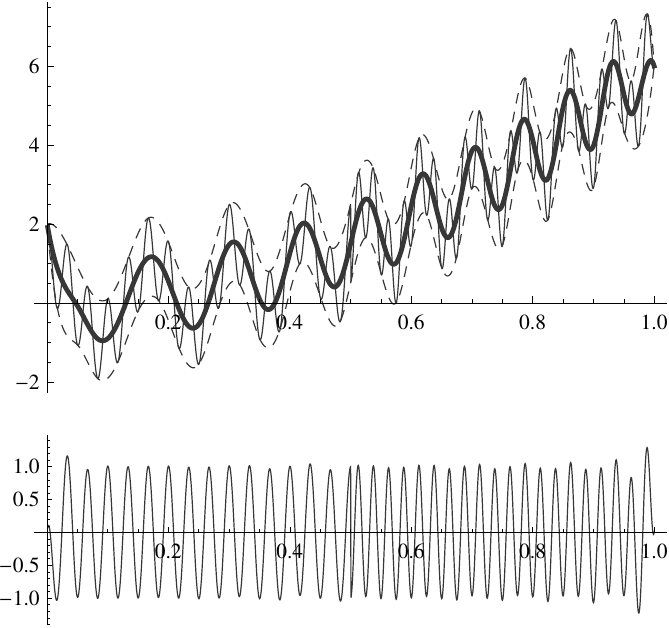}
\caption{EMD: basic IMF detection. Envelopes detection on top (thin continuous: $f$, dashed: $\bar{f}$ and $\underline{f}$ and thick continuous: $m$). 
On bottom: the first IMF candidate $r_1$.}
\label{fig:imf}
\end{figure}

In 1998, Huang et al. \cite{Huang1998} proposed an original method called Empirical Mode Decomposition (EMD) to decompose a signal into specific modes (we define the meaning of ``mode'' hereafter).
Its particularity is that it does not use any prescribed function basis but it is self adapting accordingly to the analyzed signal $f(t)$. 
In this paper, as we will use the Fourier formalism in section~\ref{sec:ewt}, we adopt the description used in \cite{Daubechies2011} which is 
slightly different from the original used in \cite{Huang1998}.
EMD aims to decompose a signal as a (finite) sum of $N+1$ \textit{Intrinsic Mode Functions} (IMF) $f_k(t)$ such that 
\begin{equation}
f(t)=\sum_{k=0}^Nf_k(t).
\end{equation}
An IMF is an amplitude modulated-frequency modulated function which can be written in the form
\begin{equation}
f_k(t)=F_k(t)\cos\left(\varphi_k(t)\right) \qquad \text{where}\; F_k(t),\varphi_k'(t)>0 \quad \forall t.
\end{equation}
The main assumption is that $F_k$ and $\varphi_k'$ vary much slower than $\varphi_k$. The IMF $f_k$ behaves as a harmonic component.
Originally, the method of Huang et al. \cite{Huang1998} to extract such IMFs is a pure algorithmic method. Candidates for an IMF are extracted by first 
computing the upper, $\bar{f}(t)$, and lower, $\underline{f}(t)$, envelopes via a cubic spline interpolation from the maxima and minima of $f$. 
Then the mean envelope is obtained by computing $m(t)=(\bar{f}(t)+\underline{f}(t))/2$ and finally the candidate by $r_1(t)=f(t)-m(t)$ (see Fig.~\ref{fig:imf}). Generally, $r_1(t)$ does not 
fulfill the properties of an IMF. A good candidate can be reached by iterating the same process to $r_1$ and the subsequent $r_k$. The final retained IMF is $f_1(t)=r_n(t)$. 
Then the next IMF is obtained by the same algorithm applied on $f(t)-f_1(t)$. The remaining IMFs can be computed by repeating this algorithm on the successive residues. 

The interesting fact about this algorithm is that it is highly adaptable and is able to extract the non-stationary part of the original function. However, its main problem is that it is based 
on an ad-hoc process which is mathematically difficult to model. Consequently it is difficult to really understand what the EMD provides. For example, some problems appear when some noise is present in the signal. 
To deal with this problem, an Ensemble EMD (EEMD) was proposed in 
\cite{Torres2011}. The authors propose to compute several EMD decompositions of the original signal corrupted by different artificial noises. Then the final EEMD is the average 
of each EMD. This approach seems to stabilize the obtained decomposition but it increases the computational cost. 

Another EMD approach is proposed in \cite{Hou2011}. The authors proposed to minimize a functional which looks for a sparse representation of 
$f$ in a dictionary of IMFs. This variational method provides similar results as the original EMD algorithm. However, this functional is based on a scheme which 
uses higher order total variation terms, this makes the method sensitive to the presence of noise and some filtering must be added to the method.

\subsection{Wavelets approaches}\label{sec:wav}
\begin{figure}[!t]
\centering
\includegraphics[width=\columnwidth]{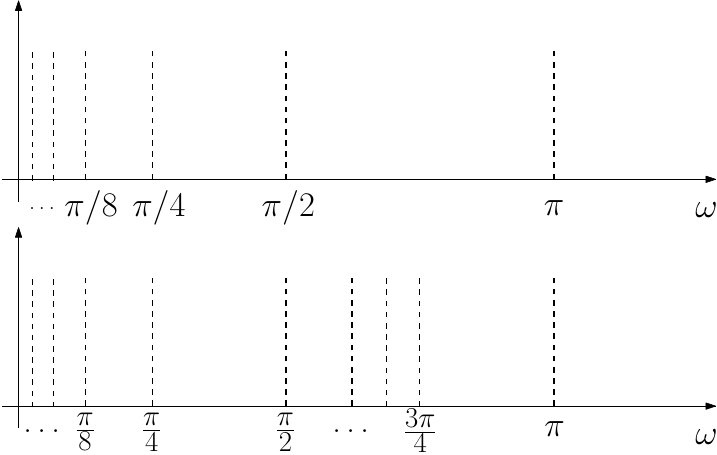}
\caption{On top: dyadic wavelet tiling of the frequency line. On bottom: a wavelet packet like tiling.}
\label{fig:dyadictiling}
\end{figure}

Nowadays, wavelet analysis is one of the most used tool in signal analysis. Let us fix some notations and recall the very basics about wavelet theory. For further details, 
we refer the reader to the extensive literature about the wavelet theory, see for example \cite{Daubechies1992,Jaffard2001,Mallat2009,Meyer1997a}. The Fourier transform and its inverse are denoted 
$\hat{f}$ and $\check{f}$, respectively. In the temporal domain, a wavelet dictionary $\{\psi_{u,s}\}$ is defined as the dilated, with a parameter $s>0$, and translated by $u\in\R$ of a mother wavelet $\psi$ 
(of zero-mean) as
\begin{equation}
\psi_{u,s}(t)=\frac{1}{\sqrt{s}}\psi\left(\frac{t-u}{s}\right).
\end{equation}
Then the wavelet transform of $f$ is obtained by computing the inner products $\mathcal{W}_f(u,s)=\langle f,\psi_{u,s}\rangle$. If $s$ is a continuous variable then $\mathcal{W}_f(u,s)$ is called the continuous wavelet transform while 
if $s=a^j$ then $\mathcal{W}_f(u,s)=\mathcal{W}_f(u,j)$ is called the discrete wavelet transform. A useful property of the wavelet transform is that it can be viewed as the application of a filter bank (each filter corresponds to one scale).
In practice, the most used case is the dyadic case, $s=2^j$. 
It can be shown that such a case corresponds to tile the time-frequency plane like in top of Fig.~\ref{fig:dyadictiling}.

As we are interested in developing adaptive representations, we recall some existing tentatives of adaptive wavelets construction. As far as we know there are a very few attempts in the literature. Probably the
most known method is the wavelet packets in a basis pursuit framework based on successive scale refinements of the expansion. It provides an adaptive time-frequency plane tiling like in bottom of Fig.~\ref{fig:dyadictiling}.
Even though the wavelet packets are useful in many applications, they use a constant prescribed ratio in the subdivision scheme, which limits their adaptability.\\
Another approach, called the Malvar-Wilson wavelets \cite{Jaffard2001,Malvar1990}, tries to build an adaptive representation by segmenting the temporal signal itself in order to separate the time intervals containing different spectral information. 
While the original idea is interesting, it turns out that the temporal segmentation is a difficult task to perform efficiently.\\
In \cite{Meyer1997}, the authors propose a method, called the brushlets, which aims to build an adaptive filter bank directly in the Fourier domain. Basically, it uses the idea of the Malvar-Wilson wavelets but 
segments the Fourier spectrum of the signal, instead of the signal itself. Conceptually the ideas in this work are really interesting, however the proposed construction is quite complicated and is 
also based on prescribed subdivisions. \\
The last work we want to mention is a recent work of Daubechies et al. \cite{Daubechies2011} entitled ``synchrosqueezed wavelets''. This approach combines a classic wavelet analysis and a reallocation method of the 
time-frequency plane information. 
This algorithm permits to obtain a more accurate time-frequency representation and consecutively it is possible to extract specific ``modes'' by choosing the appropriate information to keep.\\
All the above methods use either a prescribed scale subdivision schemes or a smart utilization of the output of a classic wavelet analysis. As far as we know, no work exists which aims to build a full adaptive wavelet transform. 
The remaining of the paper will addresses such construction.

\section{Empirical Wavelets}\label{sec:ewt}
\subsection{Definition}
We propose a method to build a family of wavelets adapted to the processed signal. If we take the Fourier point of view, this construction is equivalent to building a set of bandpass filters. One way to reach 
the adaptability is to consider that the filters' supports depend on where the information in the spectrum of the analyzed signal is located. Indeed, the IMF properties are equivalent to say that 
the spectrum of an IMF is of compact support and centered around a specific frequency (signal dependent).
For clarity, we only consider real signals (their spectrum is symmetric with respect to the frequency $\omega=0$) but the following reasoning can be easily extended to complex signal by building different filters in
the positive and negative frequencies, respectively. We also consider a normalized Fourier axis which have a $2\pi$ periodicity, in order to respect the Shannon criteria, and we restrict our discussion to $\omega\in[0,\pi]$.

\begin{figure}[!t]
\centering
\includegraphics[width=\columnwidth]{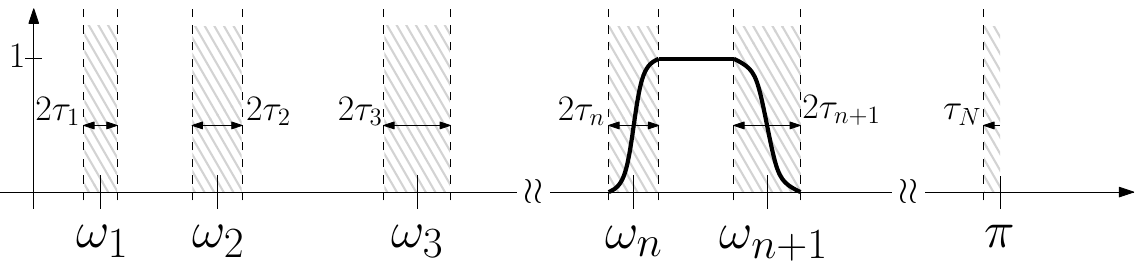}
\caption{Partitioning of the Fourier axis}
\label{fig:tiling}
\end{figure}

Let us start by assuming that the Fourier support $[0,\pi]$ is segmented into $N$ contiguous segments (we will discuss later how we can obtain such partitioning). We denote $\omega_n$ to be the limits between each segments 
(where $\omega_0=0$ and $\omega_N=\pi$), see Fig.~\ref{fig:tiling}. 
Each segment is denoted $\Lambda_n=[\omega_{n-1},\omega_n]$, then it is easy to see that $\bigcup_{n=1}^N\Lambda_n=[0,\pi]$. Centered around each $\omega_n$, we define a transition phase 
(the gray hatched areas on Fig.~\ref{fig:tiling}) $T_n$ of width $2\tau_n$. 

\begin{figure}[!t]
\centering
\includegraphics[width=\columnwidth]{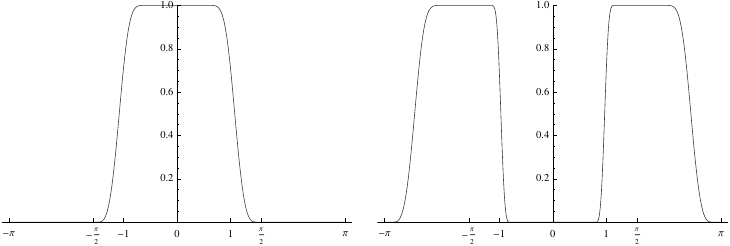}
\caption{On left: Fourier transform of the scaling function for $\nu_n=1,\gamma=0.5$. On right: Fourier transform of the wavelet function for $\nu_n=1,\nu_{n+1}=2.5,\gamma=0.2$.}
\label{fig:phipsiex}
\end{figure}

The empirical wavelets are defined as bandpass filters on each $\Lambda_n$.
To do so, we utilize the idea used in the construction of both Littlewood-Paley and Meyer's wavelets \cite{Daubechies1992}. 
Then $\forall n>0$, we define the empirical scaling function and the empirical wavelets by expressions of equations (\ref{eq:phi}) and (\ref{eq:psi}), respectively. 
\begin{equation}\label{eq:phi}
\hat{\phi}_n(\omega)=
\begin{cases}
1 \quad \qquad \qquad \text{if}\;|\omega|\leq \omega_n-\tau_n\\
\cos\left[\frac{\pi}{2}\beta\left(\frac{1}{2\tau_n}(|\omega|-\omega_n+\tau_n)\right)\right] \\
\; \; \quad \qquad \qquad  \text{if}\;\omega_n-\tau_n\leq |\omega|\leq\omega_n+\tau_n \\
0 \quad \qquad \qquad \text{otherwise}
\end{cases}
\end{equation}

and

\begin{equation}\label{eq:psi}
\hat{\psi}_n(\omega)=
\begin{cases}
1  \qquad\qquad \text{if}\; \omega_n+\tau_n\leq |\omega|\leq \omega_{n+1}-\tau_{n+1} \\
\cos\left[\frac{\pi}{2}\beta\left(\frac{1}{2\tau_{n+1}}(|\omega|-\omega_{n+1}+\tau_{n+1})\right)\right]  \\ 
\qquad\qquad\;\; \text{if}\; \omega_{n+1}-\tau_{n+1}\leq |\omega|\leq \omega_{n+1}+\tau_{n+1}\\
\sin\left[\frac{\pi}{2}\beta\left(\frac{1}{2\tau_n}(|\omega|-\omega_n+\tau_n)\right)\right] \\
\qquad\qquad\;\; \text{if}\; \omega_n-\tau_n\leq |\omega|\leq \omega_n+\tau_n\\
0 \qquad\qquad \text{otherwise.}
\end{cases}
\end{equation}
The function $\beta(x)$ is an arbitrary $\mathcal{C}^k([0,1])$ function such that
\begin{equation}\label{eq2}
\beta(x)=
\begin{cases}
0 & \text{if}\; x\leq 0\\
1 & \text{if}\; x\geq 1
\end{cases}
\; \text{and}\; \beta(x)+\beta(1-x)=1 \quad\forall x\in [0,1].
\end{equation}
Many functions satisfy these properties, the most used in the literature \cite{Daubechies1992} is 
\begin{equation}
\beta(x)=x^4(35-84x+70x^2-20x^3).
\end{equation}
Concerning the choice of $\tau_n$, several options are possible. The simplest is to choose $\tau_n$ proportional to $\omega_n$: $\tau_n=\gamma\omega_n$ where $0<\gamma<1$. 
Consequently, $\forall n>0$, Eq.~(\ref{eq:phi}) and (\ref{eq:psi}) simplify to Eq.~(\ref{eq:phi2}) and (\ref{eq:psi2})

\begin{equation}\label{eq:phi2}
\hat{\phi}_n(\omega)=
\begin{cases}
1  \qquad\qquad\qquad \text{if}\;|\omega|\leq (1-\gamma)\omega_n\\
\cos\left[\frac{\pi}{2}\beta\left(\frac{1}{2\gamma\omega_n}(|\omega|-(1-\gamma)\omega_n)\right)\right] \\
\qquad\qquad\qquad\;\; \text{if}\; (1-\gamma)\omega_n\leq |\omega|\leq (1+\gamma)\omega_n \\
0  \qquad\qquad\qquad \text{otherwise,}
\end{cases}
\end{equation}

and

\begin{equation}\label{eq:psi2}
\hat{\psi}_n(\omega)=
\begin{cases}
1 \qquad\qquad\;\; \text{if}\; (1+\gamma)\omega_n\leq |\omega|\leq (1-\gamma)\omega_{n+1} \\
\cos\left[\frac{\pi}{2}\beta\left(\frac{1}{2\gamma\omega_{n+1}}(|\omega|-(1-\gamma)\omega_{n+1})\right)\right]  \\ 
\qquad\qquad\quad \text{if}\; (1-\gamma)\omega_{n+1}\leq |\omega|\leq (1+\gamma)\omega_{n+1}\\
\sin\left[\frac{\pi}{2}\beta\left(\frac{1}{2\gamma\omega_n}(|\omega|-(1-\gamma)\omega_n)\right)\right] \\
 \qquad\qquad\quad\text{if}\; (1-\gamma)\omega_n\leq |\omega|\leq (1+\gamma)\omega_n\\
0 \qquad\qquad\;\, \text{otherwise.}
\end{cases}
\end{equation}
An example of $\hat{\phi}_n$ for $\nu_n=1,\gamma=0.5$ and $\hat{\psi}_n$ for $\nu_n=1,\nu_{n+1}=2.5,\gamma=0.2$ is given in Fig.~\ref{fig:phipsiex}.

\subsection{Segmentation of the Fourier spectrum}
How we segment the Fourier spectrum is important as it is the step which provides the adaptability with respect to the analyzed signal to our method. We aim to separate different portions 
of the spectrum which correspond to modes e.g centered around a specific frequency and of compact support. In this paper, we assume that the number of segment, $N$, is given (in section~\ref{sec:detect} we 
propose a method to estimate the number of bands). This implies that we need a total of 
$N+1$ boundaries, but 0 an $\pi$ are always used in our definition and consequently we need to find $N-1$ extra boundaries.
To find such boundaries, we first detect the local maxima in the spectrum and sort them in decreasing order (0 and $\pi$ are excluded). Let us assume that the algorithm found $M$ maxima. Two cases can appear:
\begin{itemize}
 \item $M\geq N$: the algorithm found enough maxima to define the wanted number of segments, then we keep only the first $N-1$ maxima,
 \item $M<N$: the signal has less modes than expected, then we keep all the detected maxima and reset $N$ to the appropriate value.
\end{itemize}
Now, equipped with this set of maxima plus 0 and $\pi$, we define the boundaries $\omega_n$ of each segment as the center between two consecutive maxima. 

\subsection{Frame}
The following proposition shows that, by properly choosing the parameter $\gamma$, we can obtain a tight frame.

\begin{proposition}
If $\gamma<\min_n\left(\frac{\omega_{n+1}-\omega_n}{\omega_{n+1}+\omega_n}\right)$, then the set $\{\phi_1(t),\{\psi_n(t)\}_{n=1}^N\}$ is a tight frame of $L^2(\R)$.
\end{proposition}

\begin{IEEEproof}
We follow the idea behind the construction of Meyer's wavelet.\\
The set $\{\phi_1(t),\{\psi_n(t)\}_{n=1}^N\}$ is a tight frame if
\begin{equation}
\sum_{k=-\infty}^{+\infty}\left(\left|\hat{\phi_1}(\omega+2k\pi)\right|^2+\sum_{n=1}^N\left|\hat{\psi}_n(\omega+2k\pi)\right|^2\right)=1.
\end{equation}

\begin{figure}[!t]
\centering
\includegraphics[width=\columnwidth]{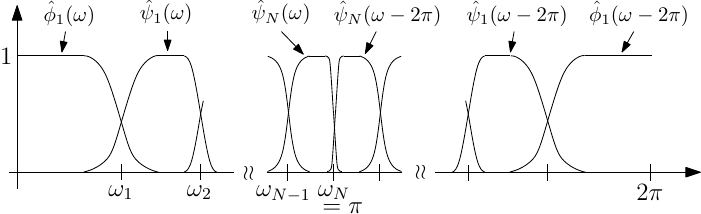}
\caption{Periodicity of the filter bank}
\label{fig:perio}
\end{figure}

Accordingly to the $2\pi$ periodicity (see Fig.~\ref{fig:perio}), it is enough to focus on the interval $[0,2\pi]$. Following the previous notations, we can write
\begin{equation}
[0,2\pi]=\bigcup_{n=1}^N\Lambda_n\cup\bigcup_{n=1}^N\Lambda_{\sigma(n)},
\end{equation}
where $\Lambda_{\sigma(n)}$ is a copy of $\Lambda_n$ but centered at $2\pi-\nu_n$ instead $\nu_n$.
First, it is easy to see that for 
$\omega\in\left(\bigcup_{n+1}^N\Lambda_n/\bigcup_{n+1}^NT_n\right)\cup\left(\bigcup_{n=1}^N\Lambda_{\sigma(n)}/\bigcup_{n+1}^NT_{\sigma(n)}\right)$ we have
\begin{align}
\left|\hat{\phi_1}(\omega)\right|^2 +&\left|\hat{\phi_1}(\omega-2\pi)\right|^2 + \\ \notag
 \sum_{n=1}^N & \left(\left|\hat{\psi}_n(\omega)\right|^2 +\left|\hat{\psi}_n(\omega-2\pi)\right|^2\right)=1.
\end{align}
Then, it remains to look at the transition areas. Because of properties of $\beta$, this result also holds in $T_n$ if consecutive $T_n$ do not overlap:
\begin{gather}
\tau_n+\tau_{n+1}<\omega_{n+1}-\omega_n\\
\Leftrightarrow \gamma\omega_n+\gamma\omega_{n+1}<\omega_{n+1}-\omega_n\\
\Leftrightarrow \gamma<\frac{\omega_{n+1}-\omega_n}{\omega_{n+1}+\omega_n} \label{eq:cond}.
\end{gather}
This condition must be true for all $n$ which is equivalent to said that condition (\ref{eq:cond}) must be true for the smallest $T_n$ and finally we get the result if 
$\gamma<\min_n\left(\frac{\omega_{n+1}-\omega_n}{\omega_{n+1}+\omega_n}\right)$.
\end{IEEEproof}

\begin{figure}[!t]
\centering
\includegraphics[width=\columnwidth,height=1in]{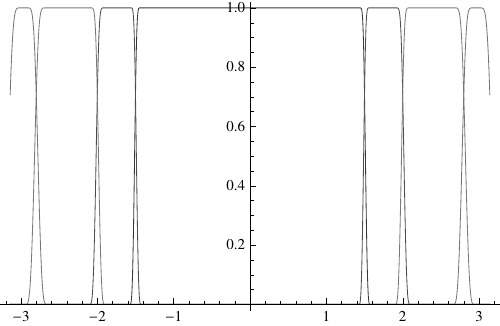} 
\caption{Example of a Fourier partitioning of an empirical filter bank (see text for details)}
\label{fig:filterbank}
\end{figure}

Fig.~\ref{fig:filterbank} gives an empirical filter bank example based on the set $\omega_n\in\{0,1.5,2,2.8,\pi\}$ with $\gamma=0.05$ (the theory tells us that $\gamma<0.057$).

\subsection{Empirical wavelet transform}
From the previous section, we know how to build a tight frame set of empirical wavelets. We can now define the Empirical Wavelet Transform (EWT), $\mathcal{W}_f^{\mathcal{E}}(n,t)$, 
in the same way as for the classic wavelet transform. The detail coefficients are 
given by the inner products with the empirical wavelets:
\begin{align}
\mathcal{W}_f^{\mathcal{E}}(n,t)=\langle f,\psi_n\rangle&=\int f(\tau)\overline{\psi_n(\tau-t)}d\tau \\
&=\left(\hat{f}(\omega)\overline{\hat{\psi}_n(\omega)}\right)^{\vee},
\end{align}
and the approximation coefficients (we adopt the convention $\mathcal{W}_f^{\mathcal{E}}(0,t)$ to denote them) by the inner product with the scaling function:
\begin{align}
\mathcal{W}_f^{\mathcal{E}}(0,t)=\langle f,\phi_1\rangle&=\int f(\tau)\overline{\phi_1(\tau-t)}d\tau \\
&=\left(\hat{f}(\omega)\overline{\hat{\phi}_1(\omega)}\right)^{\vee},
\end{align}
where $\hat{\psi}_n(\omega)$ and $\hat{\phi}_1(\omega)$ are defined by Eq.~\ref{eq:psi2} and \ref{eq:phi2}, respectively. The reconstruction is obtained by
\begin{align}
f(t)&=\mathcal{W}_f^{\mathcal{E}}(0,t)\star\phi_1(t)+\sum_{n=1}^N\mathcal{W}_f^{\mathcal{E}}(n,t)\star\psi_n(t)\\
&=\left(\widehat{\mathcal{W}_f^{\mathcal{E}}}(0,\omega)\hat{\phi}_1(\omega)+\sum_{n=1}^N\widehat{\mathcal{W}_f^{\mathcal{E}}}(n,\omega)\hat{\psi}_n(\omega)\right)^{\vee}.
\end{align}
Following this formalism, the empirical mode $f_k$, as defined in section~\ref{sec:emd}, is given by
\begin{gather}
f_0(t)=\mathcal{W}_f^{\mathcal{E}}(0,t)\star\phi_1(t),\\
f_k(t)=\mathcal{W}_f^{\mathcal{E}}(k,t)\star\psi_k(t).
\end{gather}

\section{Experiments}\label{sec:expe}
\subsection{Test signals}
We propose to test the Empirical Wavelet Transform on four different signals: three are artificial signals taken from \cite{Hou2011} and one real electrocardiogram signal. Their description is given hereafter.

\begin{figure}[!t]
\centering
\includegraphics[width=\columnwidth,height=0.7\textheight]{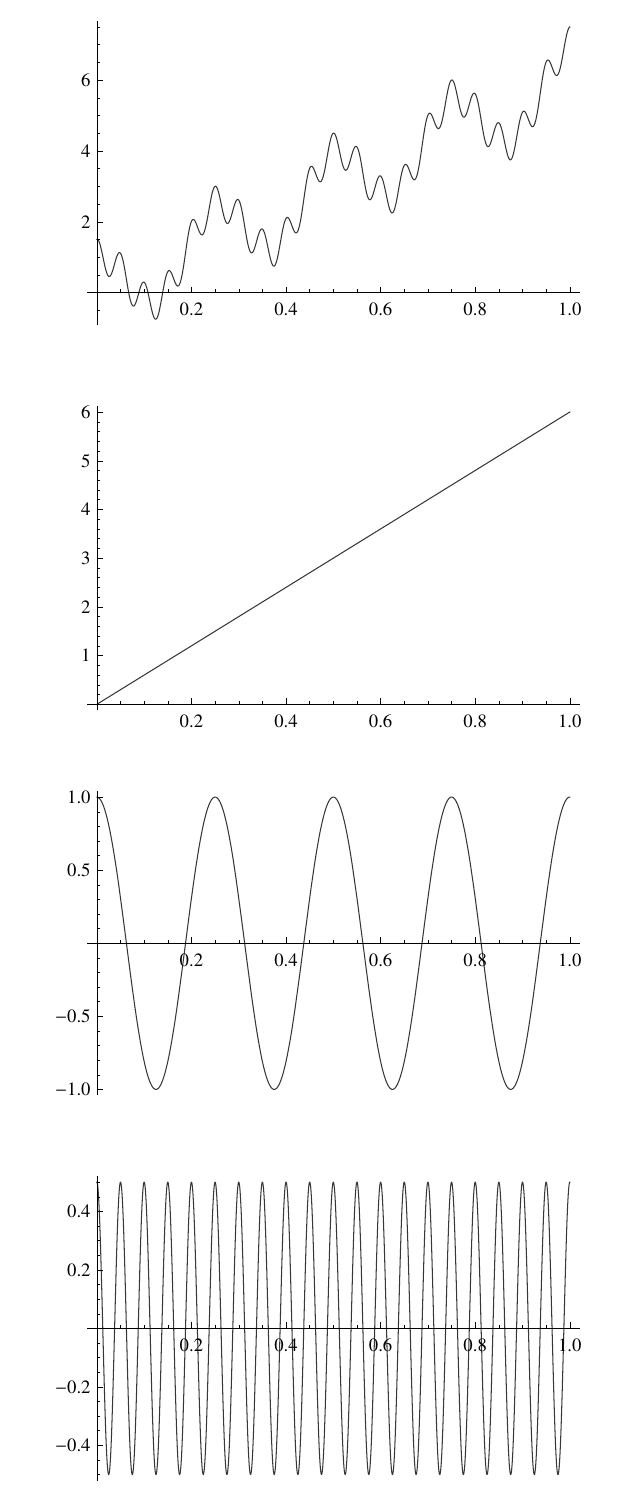}
\caption{$f_{Sig1}$ test signal: signal on top. From second to last row: the different components constituting {\it Sig1}.}
\label{fig:sig1}
\end{figure}

\paragraph{Simulated $f_{Sig1}$}
The first test signal, {\it Sig1}, is made with the sum of three distinct components (for $t\in[0,1]$) (see Fig.~\ref{fig:sig1}):
\begin{gather}
f_{c1}(t)=6t \\
f_{c2}(t)=\cos(8\pi t)\\
f_{c3}(t)=0.5\cos(40\pi t)
\end{gather}
and then
\begin{equation}
f_{Sig1}(t)=f_{c1}(t)+f_{c2}(t)+f_{c3}(t)
\end{equation}

\begin{figure}[htbp]
\centering
\includegraphics[width=\columnwidth,height=0.7\textheight]{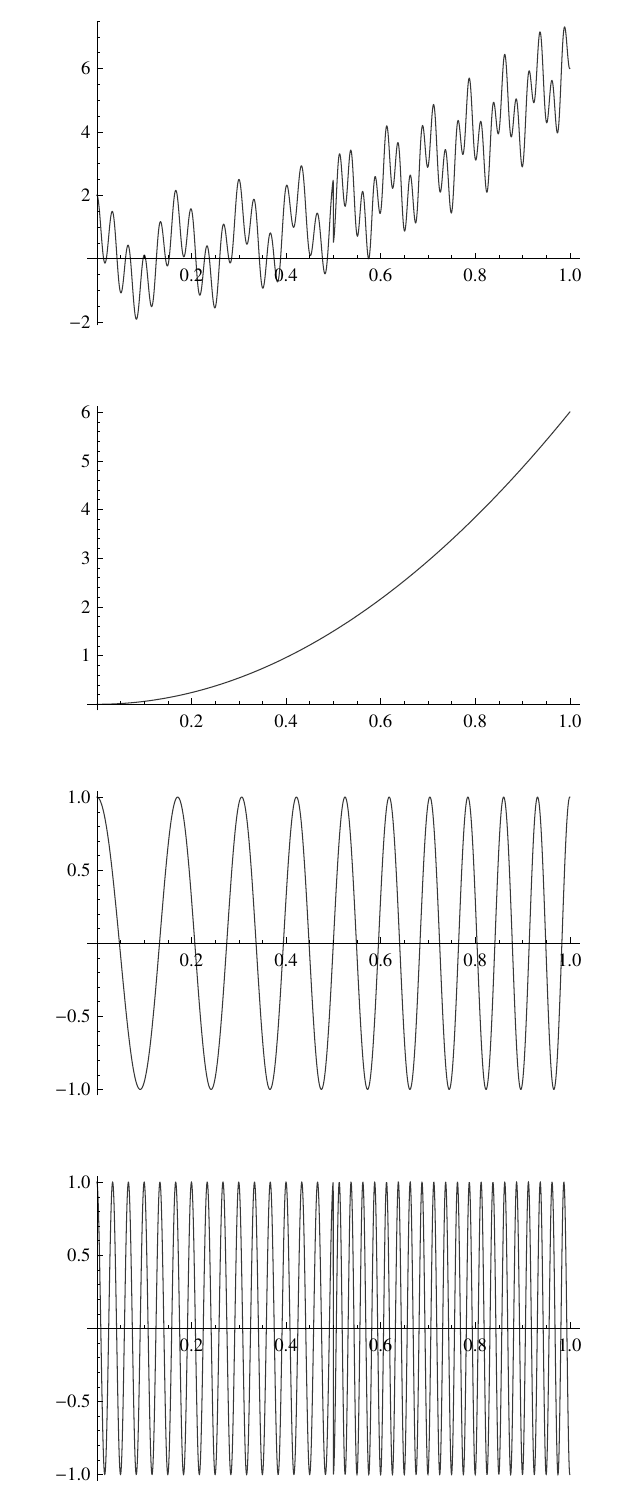}
\caption{$f_{Sig2}$ test signal: signal on top. From second to last row: the different components constituting {\it Sig2}.}
\label{fig:sig2}
\end{figure}

\paragraph{Simulated $f_{Sig2}$}
The second test signal, {\it Sig2}, is made with the sum of three distinct components (for $t\in[0,1]$) (see Fig.~\ref{fig:sig2}):
\begin{gather}
f_{c1}(t)=6t^2 \\
f_{c2}(t)=\cos(10\pi t+10\pi t^2)\\
f_{c3}(t)=
\begin{cases}
\cos(80\pi t-15\pi) & \text{if}\; t>0.5\\
\cos(60\pi t) & \text{otherwise}
\end{cases}
\end{gather}
and then
\begin{equation}
f_{Sig2}(t)=f_{c1}(t)+f_{c2}(t)+f_{c3}(t)
\end{equation}

\begin{figure}[!t]
\centering
\includegraphics[width=\columnwidth,height=0.9\textheight]{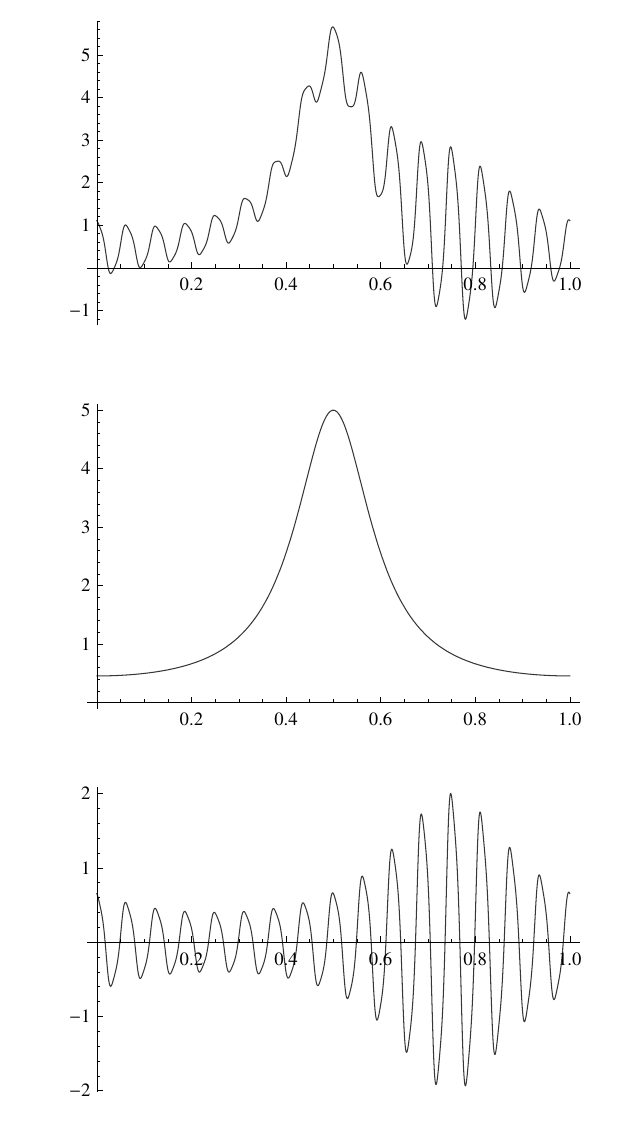}
\caption{$f_{Sig3}$ test signal: signal on top. From second to last row: $f_{c1}$ and $f_{c2}f_{c3}$, respectively.}
\label{fig:sig3}
\end{figure}

\paragraph{Simulated $f_{Sig3}$}
The third test signal, {\it Sig3}, is made with three distinct components (for $t\in[0,1]$) (see Fig.~\ref{fig:sig3}):
\begin{gather}
f_{c1}(t)=\frac{1}{1.2+\cos(2\pi t)} \\
f_{c2}(t)=\frac{1}{1.5+\sin(2\pi t)}\\
f_{c3}(t)=\cos(32\pi t+\cos(64\pi t))
\end{gather}
and then
\begin{equation}
f_{Sig3}(t)=f_{c1}(t)+f_{c2}(t)f_{c3}(t)
\end{equation}
Let notice that for this signal there is only two additive components: $f_{c1}$ and the product $f_{c2}f_{c3}$.

\begin{figure}[!t]
\centering
\includegraphics[width=\columnwidth,height=1in]{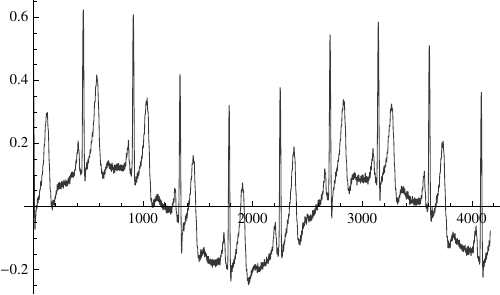}
\caption{$f_{Sig4}$: real ECG signal.}
\label{fig:sig4}
\end{figure}

\begin{figure}[!t]
\centering
\includegraphics[width=\columnwidth,height=1in]{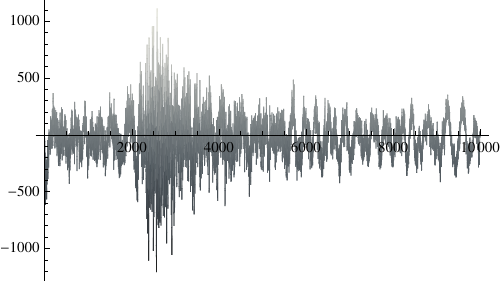}
\caption{$f_{Sig5}$: real seismic waveform signal.}
\label{fig:sig5}
\end{figure}

\paragraph{Real signals $f_{Sig4}$ and $f_{Sig5}$}
The two last signals are a real electrocardiogram (ECG) signal $f_{Sig4}$ and seismic waveform $f_{Sig5}$ illustrated in Fig.~\ref{fig:sig4} 
and Fig.~\ref{fig:sig5}, respectively.

\subsection{Comparison EMD vs. EWT}
\begin{figure}[!t]
\includegraphics[width=\columnwidth]{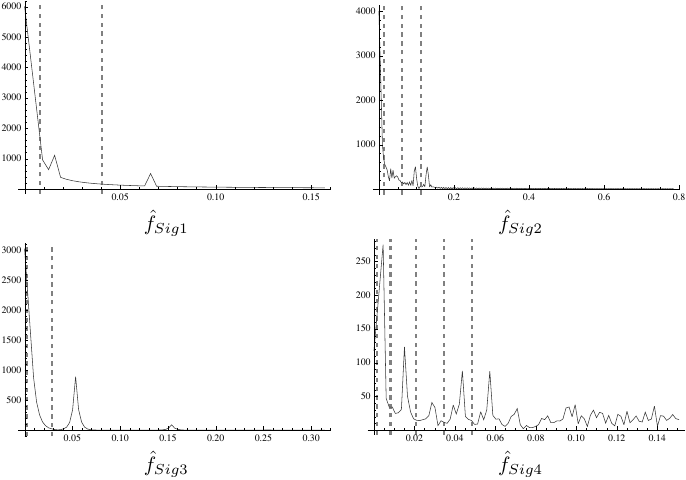}
\caption{Detected Fourier supports for signals $f_{Sig1}$ to $f_{Sig4}$.}
\label{fig:boundsig4}
\end{figure}

\begin{figure}[!t]
\centering
\includegraphics[width=\columnwidth]{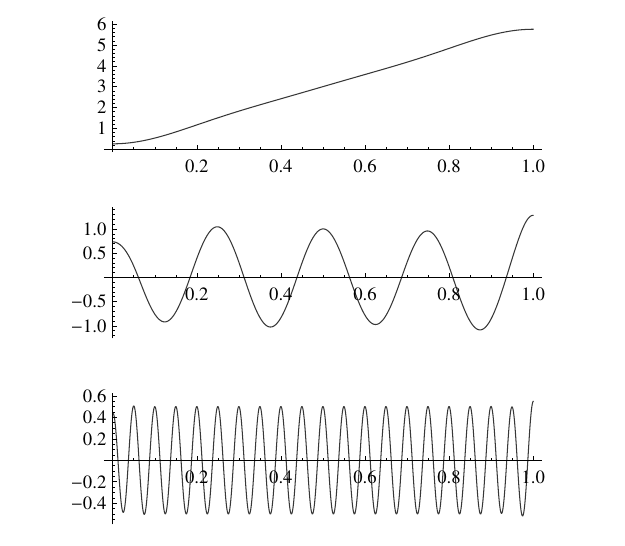}
\caption{Modes extracted by the Empirical Wavelet Transform for $f_{Sig1}$.}
\label{fig:ewtsig1}
\end{figure}

\begin{figure}[!t]
\centering
\includegraphics[width=\columnwidth]{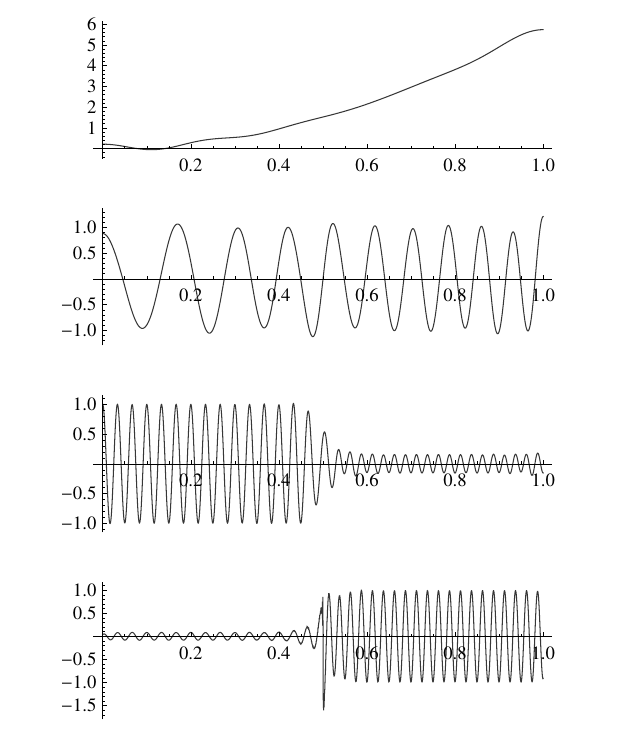}
\caption{Modes extracted by the Empirical Wavelet Transform for $f_{Sig2}$.}
\label{fig:ewtsig2}
\end{figure}

\begin{figure}[!t]
\centering
\includegraphics[width=\columnwidth,height=0.5\textheight]{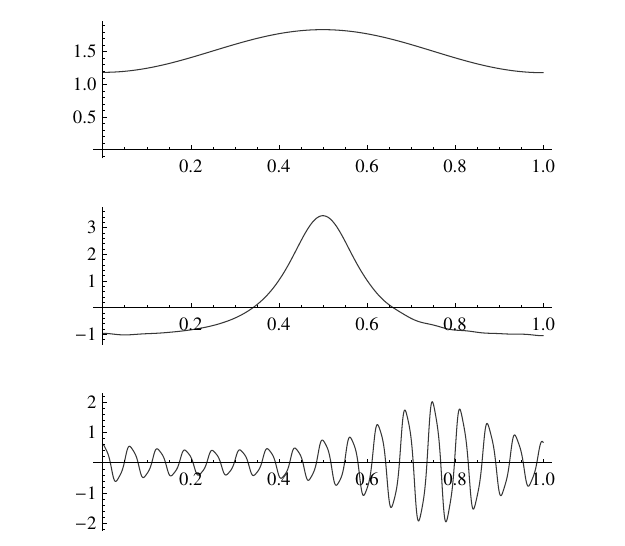}
\caption{Modes extracted by the Empirical Wavelet Transform for $f_{Sig3}$.}
\label{fig:ewtsig3}
\end{figure}

\begin{figure}[!t]
\centering
\includegraphics[width=\columnwidth,height=0.3\textheight]{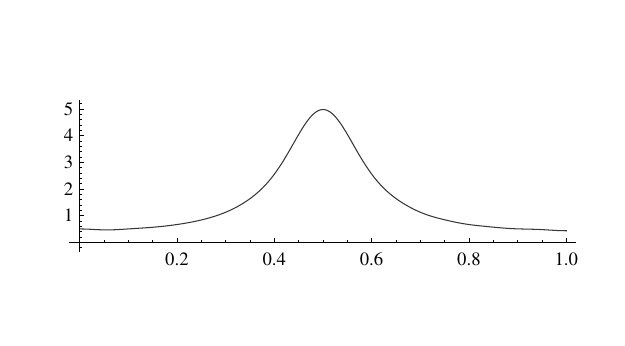}
\caption{Sum of the two first EWT modes of $f_{Sig3}$.}
\label{fig:addmode}
\end{figure}

\begin{figure}[!t]
\centering
\includegraphics[width=\columnwidth,height=0.9\textheight]{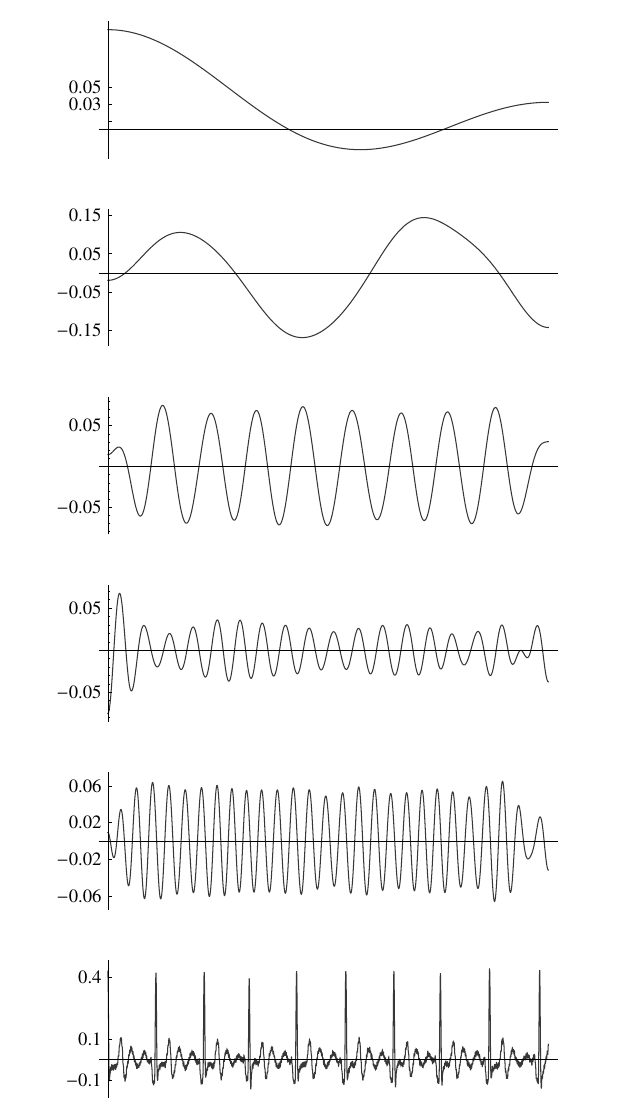}
\caption{Modes extracted by the Empirical Wavelet Transform for $f_{Sig4}$.}
\label{fig:ewtsig4}
\end{figure}

\begin{figure}[!t]
\centering
\includegraphics[width=\columnwidth,height=0.9\textheight]{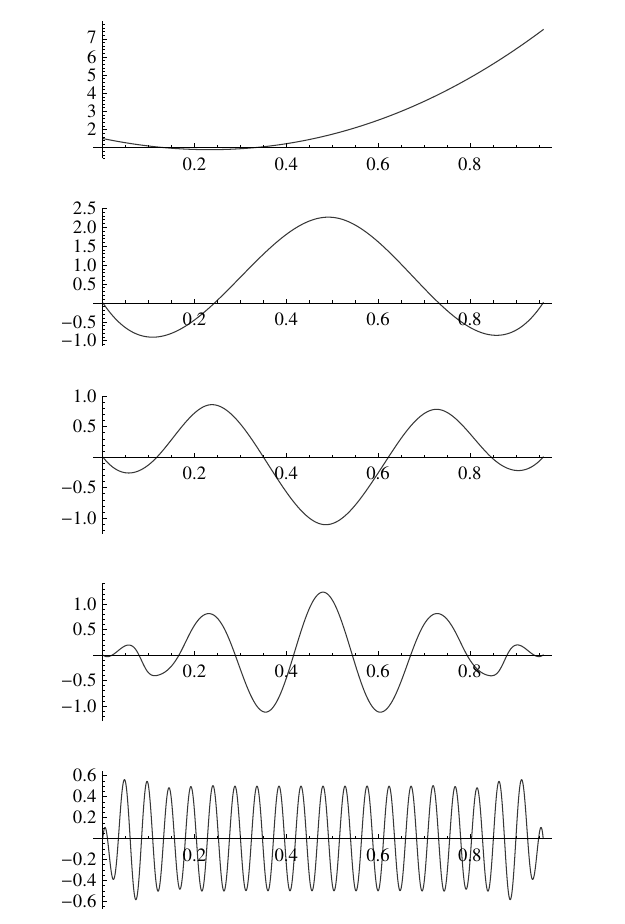}
\caption{Empirical Mode Decomposition of $f_{Sig1}$.}
\label{fig:emdsig1}
\end{figure}

\begin{figure}[!t]
\centering
\includegraphics[width=\columnwidth,height=0.9\textheight]{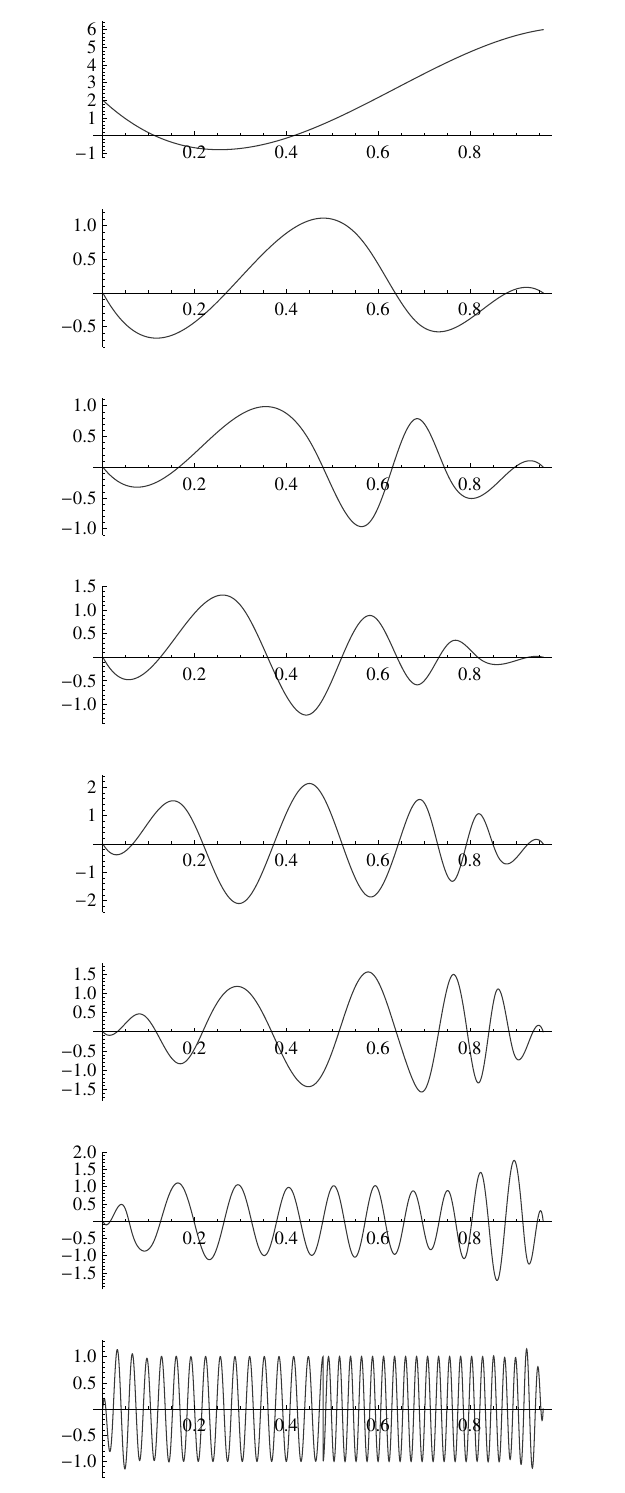}
\caption{Empirical Mode Decomposition of $f_{Sig2}$.}
\label{fig:emdsig2}
\end{figure}

\begin{figure}[!t]
\centering
\includegraphics[width=\columnwidth,height=0.9\textheight]{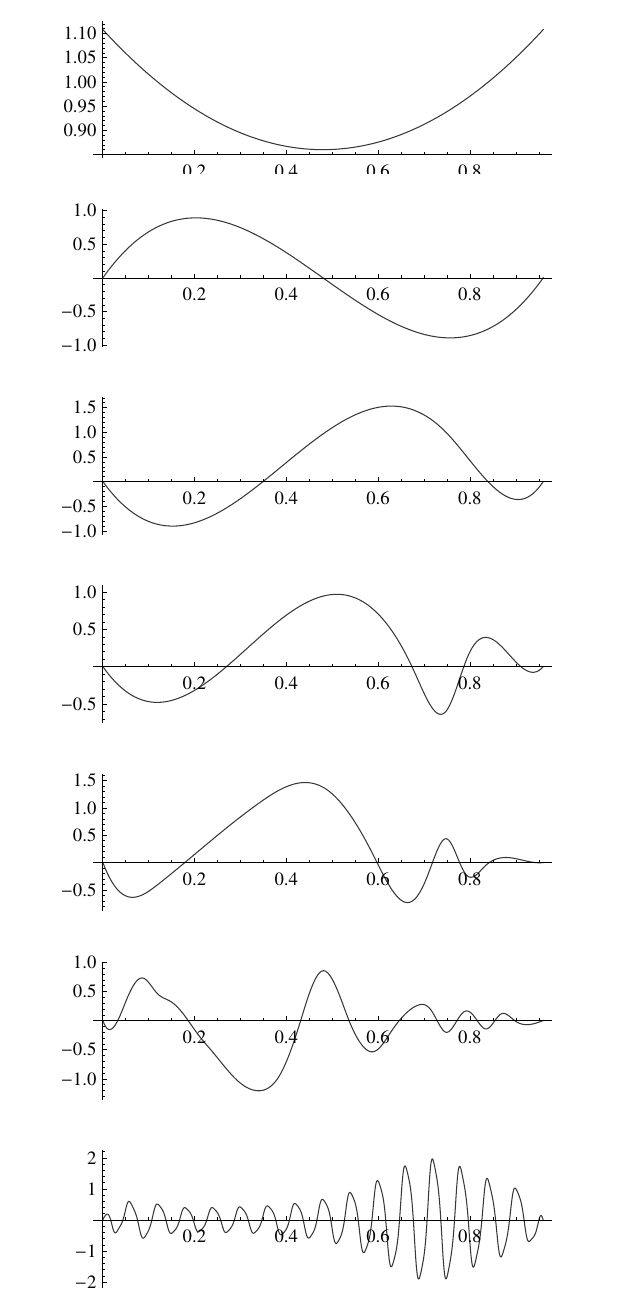}
\caption{Empirical Mode Decomposition of $f_{Sig3}$.}
\label{fig:emdsig3}
\end{figure}

\begin{figure}[!t]
\centering
\includegraphics[width=\columnwidth,height=0.9\textheight]{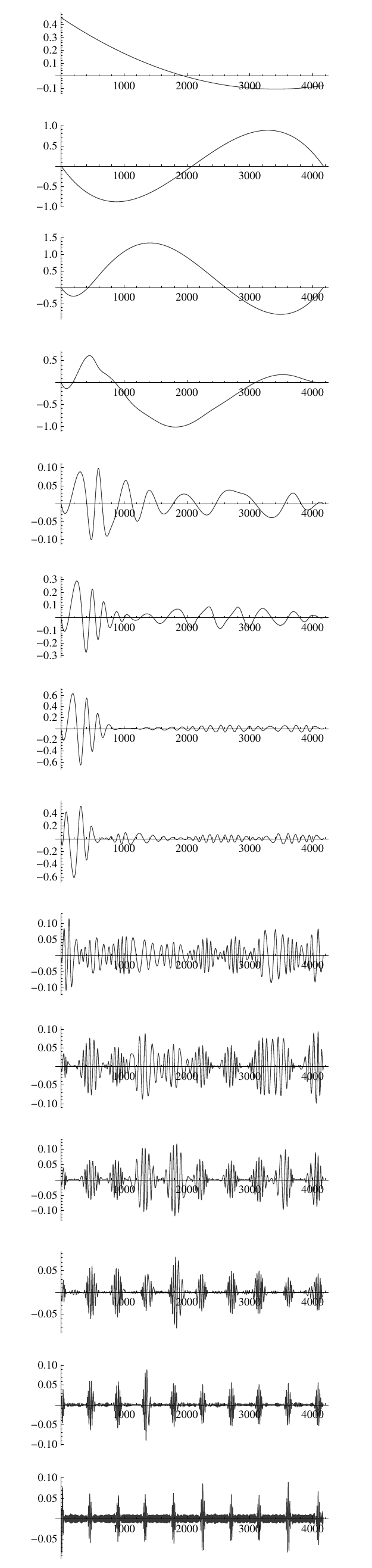}
\caption{Empirical Mode Decomposition of $f_{Sig4}$.}
\label{fig:emdsig4}
\end{figure}

In this section, we compute the Empirical Mode Decomposition and the Empirical Wavelet Transform of the four signals described in the previous section. While the EMD automatically estimate the number of modes, we 
fix a priori the number of modes, $N$, for the EWT. Consequently, the EWT output is composed of the filtering with the scaling function and the $N$ wavelets. For the different signals, we use $N_{Sig1}=2$, $N_{Sig2}=3$, 
$N_{Sig3}=2$, $N_{Sig4}=5$ and $N_{Sig5}=50$. Figures~\ref{fig:boundsig4} shows the respective spectra of each test signal and the detected boundaries for each filter support. 
We can observe that the algorithm is able is isolate the different modes. The corresponding filtered signals $\mathcal{W}_f^{\mathcal{E}}(n,t)$ are presented in Fig.~\ref{fig:ewtsig1}, \ref{fig:ewtsig2}, \ref{fig:ewtsig3} 
and \ref{fig:ewtsig4}; while the EMD outputs are given in Fig.~\ref{fig:emdsig1}, \ref{fig:emdsig2}, \ref{fig:emdsig3} and \ref{fig:emdsig4}, respectively.\\
We can observe that for the simulated test signals, the EMD always overestimates the number of modes and then separate some information which is originally part of the same component. Except for the high 
frequencies, it is difficult to interpret the EMD outputs compared with the known ``true'' components constituting the test signals. Concerning the results given by the EWT, we can see that it is able to detect 
the presence of modes in the spectrum and provides different components which are close to the original ones. However, in the case of $f_{Sig2}$, we can note that the algorithm separates the two last modes which were 
initially parts of the same component. In fact, this is not completely surprising as those modes have significant individual energy and can be considered as independent modes. About $f_{Sig3}$, another phenomenon appears. The 
initial part $f_{c1}$ is decomposed as the sum of two modes (this sum is given in Fig.~\ref{fig:addmode}). If we look more closely, the detected boundaries, shown in Fig.~\ref{fig:boundsig4}, we can see that the scaling 
function is built upon a very small support while the next support contains the main information of $f_{c1}$. If we ask the EWT to consider only $N=1$ (as it supposed to be from the construction of $f_{Sig3}$), the first 
two supports shown in Fig.~\ref{fig:boundsig4} are merged except that its boundary is moved closer to the zero frequency. This change in the boundary position has the consequence that some information of $f_{c1}$ is now 
on the second support and then appear in the second mode instead the first one. In some sense, if we ask to use $N=1$, we get a perturbed decomposition. This issue comes from the method we use to detect the boundaries 
of the Fourier support and suggests more investigations.\\
Experiments on the real ECG signal seem to give the advantage to the EWT because the EMD provides too many modes. Typically, the EMD modes six to nine are really difficult to interpret as such behavior is clearly not visible 
in the signal itself. A contrary, the EWT focuses on the oscillating patterns we can observe in {\it $f_{Sig4}$}. Of course it will be very interesting to have the opinion of a cardiologist about the medical interpretation of 
such components. The different modes for the seismic signal are not provided here but we will present its time-frequency representation in section~\ref{sec:tfr} as it 
more relevant for analysis purposes.


\section{Time-frequency representation}\label{sec:tfr}
\begin{figure}[!t]
\begin{center}
\includegraphics[width=\columnwidth]{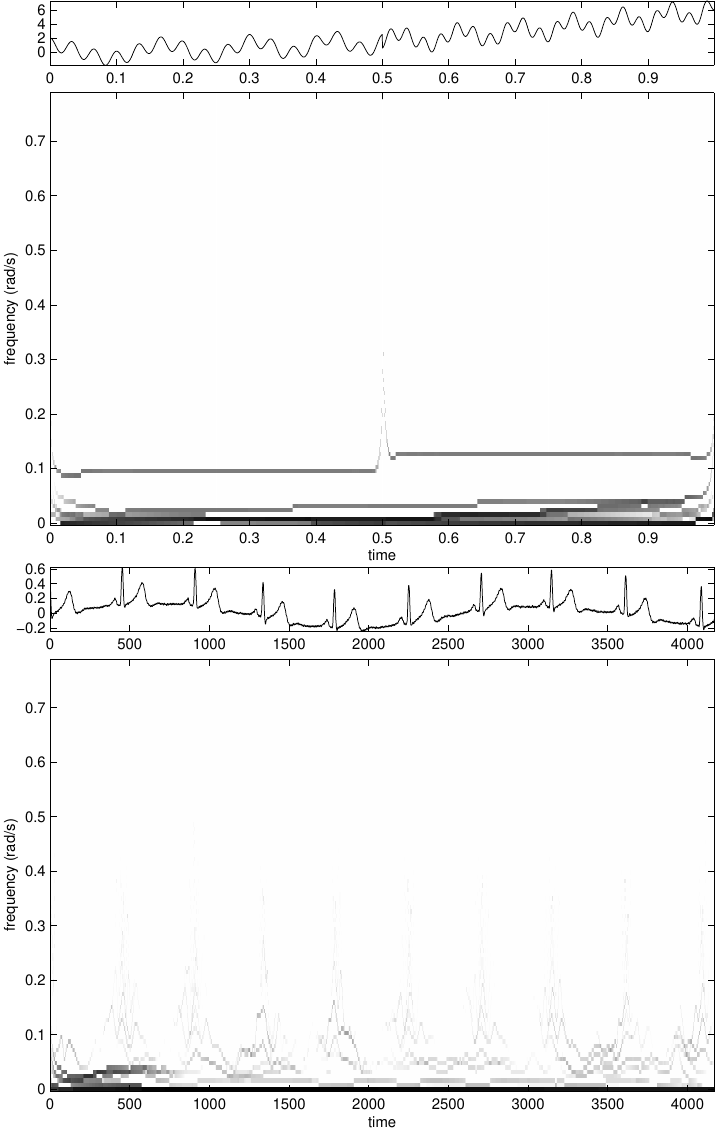} 
\end{center}
\caption{Time-frequency representation based on the EMD (Hilbert-Huang transform) for signals $f_{sig2}$ (top) and $f_{sig4}$ (bottom).}
\label{fig:TFEMD}
\end{figure}

\begin{figure}[!h]
\begin{center}
\includegraphics[width=\columnwidth]{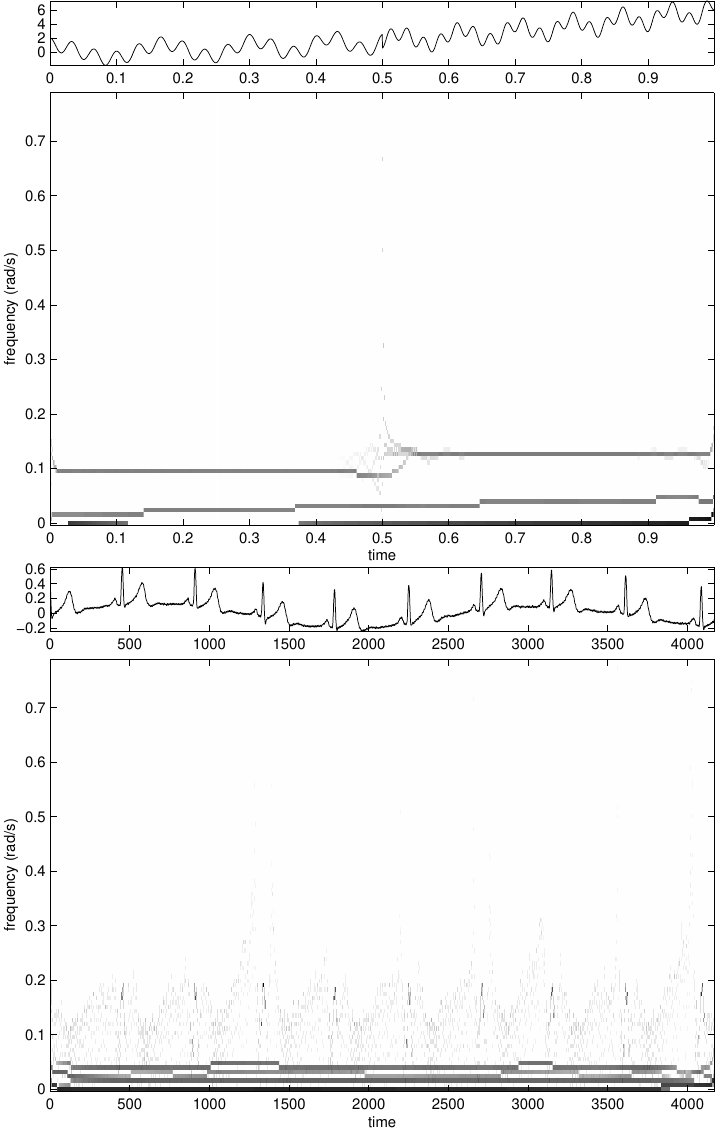} 
\end{center}
\caption{Time-frequency representation based on the EWT for signals $f_{sig2}$ (top) and $f_{sig4}$ (bottom).}
\label{fig:TFEWT}
\end{figure}

\begin{figure}[!h]
\begin{center}
\includegraphics[width=\columnwidth]{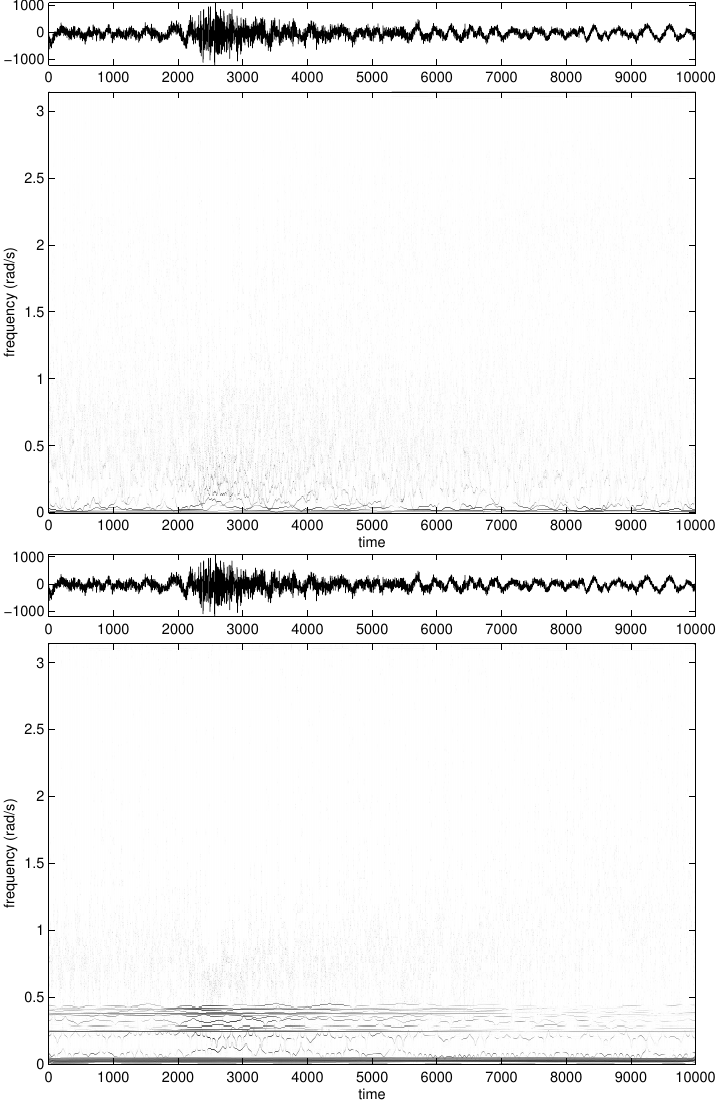} 
\end{center}
\caption{Time-frequency representation based on the EMD (on top) and the EWT (on bottom) for the seismic signal.}
\label{fig:TFPSeismic}
\end{figure}

The time-frequency representation is useful to have the information of all components summarized in a single domain. In this paper, we follow the idea used in the Hilbert-Huang transform \cite{Huang1998}. First, let us 
recall the definition of the Hilbert transform of a function $f$: 
\begin{equation}
\mathcal{H}_f(t)=\frac{1}{\pi}p.v.\int_{-\infty}^{+\infty}\frac{f(\tau)}{t-\tau}d\tau
\end{equation}
where the integral is defined by using the Cauchy principal value ($p.v.$) \cite{King}. The Hilbert transform can be used to derive the analytical form $f_a$ of $f$: 
$f_a(t)=f(t)+\imath\mathcal{H}_f(t)$. The Hilbert transform has the interesting property that in the case of AM-FM signals $f(t)=F(t)\cos\varphi(t)$, it provides 
$f_a(t)=F(t)e^{\imath\varphi(t)}$ which permits to extract the instantaneous amplitude $F(t)$ and frequency $\varphi'(t)$. The Hilbert-Huang transform 
consists in applying the Hilbert transform to each extracted IMF and then plotting each curve $\varphi_k'(t)$ in the time-frequency plane where the intensity 
of the plot is given by $F_k(t)$. For a time-frequency representation based on the EWT, we can apply the same process by evaluating the Hilbert transform 
of each filter output. For instance, figures~\ref{fig:TFEMD} and \ref{fig:TFEWT} give the time-frequency representation of $f_{sig2}$ and $f_{sig4}$ based 
on the EMD and EWT, respectively. We can notice that in the case of the Hilbert-Huang transform some artifacts appear in the low frequencies while 
representations based on the EWT avoid such problems. These artifacts are mainly due to the fact that the EMD decomposition forces the extraction of 
IMFs even if the initial components are not. Figure~\ref{fig:TFPSeismic} gives the time-frequency representations obtained from the EMD and the EWT of the 
seismic signal $f_{sig5}$ given in Fig.~\ref{fig:sig5}. For such signals the EWT seems to provide a more consistent representation during the earthquake while the 
EMD shows a sparse representation which is difficult to interpret.

\section{Automatic detection of the number of modes}\label{sec:detect}
If for a given class of signals it should be possible to guess the best number of mode $N$, this is not the case in general where no a priori information are available. In such cases, it should be interesting to 
estimate the appropriate number of modes. In general such task is difficult, hereafter we present a simple way to estimate $N$ but a deeper analysis is necessary to provide a more robust method. We 
follow the idea that the most important maxima in the magnitude of the Fourier transform of the input signal (corresponding to the center of each desired Fourier segments) are significantly larger than the 
other existing maxima. Let denote $\{M_i\}_{k=1}^M$ the set of $M$ detected maxima in the magnitude of the Fourier spectrum. Without any loss of generality we assume that this set is sorted by decreasing values 
($M_1\geq M_2\geq\ldots M_M$) and normalized in $[0;1]$. In this framework, the above idea is equivalent to keep all maxima which are greater than some amount of the difference between the bigger maximum and the smaller 
maximum. We can formulate it as ``keeping all maxima larger than the threshold $M_M+\alpha(M_1-M_M)$'' where $\alpha$ corresponds to the relative amplitude ratio. Table~\ref{tab1} gives the number of detected bands for 
each signal according to different values of $\alpha$. We can see that values of $\alpha$ around 0.3 and 0.4 seems to give consistent results which correspond to a trade off between too much detection and a good separation of 
the information in the Fourier spectrum.

\begin{table}
\begin{center}
\begin{tabular}{|c|c|c|c|c|} \hline
$\alpha$ & $f_{sig1}$ & $f_{sig2}$ & $f_{sig3}$ & $f_{sig4}$ \\ \hline
0.1 & 4 & 17 & 7 & 24 \\ \hline
0.3 & 3 & 9 & 4 & 5 \\ \hline
0.4 & 3 & 9 & 4 & 4 \\ \hline
0.5 & 2 & 7 & 4 & 3 \\ \hline
\end{tabular}
\end{center}
\caption{Number of bands automatically detected for different $\alpha$ for each signal.}
\label{tab1}
\end{table}

\section{Extension to images}\label{sec:2d}
In this section, we address the opportunity to extend the previous EWT to 2D signals (images). The basic idea to built such 2D-EWT is to use a tensor product 
as for classic wavelets, which means process the rows and then the columns of the input image by using the 1D-EWT defined previously. If this idea seems obvious 
to adopt, it turns out that two main issues arise: there is no guarantee that two different rows (or columns) will have the same number of Fourier supports; 
even if we have the same number of bands for each rows (or columns), the supports corresponding to the same band number can be very far away from each other. 
The direct consequence is that in the output image corresponding to a specific band, very different frequency modes can be mixed providing a ``strange'' 
representation. In order to avoid these issues, a solution consists to use the same filters for each rows (or columns). A simple way to build such solution 
is to consider a mean spectrum for the rows (or the columns), then perform the detection of the Fourier supports based on this mean spectrum and finally 
use the same filters for all rows (or columns). The complete algorithm corresponds to perform the following steps (we denote $f$ the input image):
\begin{enumerate}
\item Take the 1D FFT of each rows $i$ of $f$; $\hat{f}(i,\omega)$; and compute the mean row spectrum magnitude: $\widetilde{\mathcal{F}}_{row}=\frac{1}{N_{row}}\sum_{i=0}^{N_{row}}\hat{f}(i,\omega)$,
\item Take the 1D FFT of each columns $j$ of $f$; $\hat{f}(\omega,j)$; and compute the mean row spectrum magnitude: $\widetilde{\mathcal{F}}_{columns}=\frac{1}{N_{columns}}\sum_{j=0}^{N_{columns}}\hat{f}(\omega,j)$,
\item Perform the boundaries detection on $\widetilde{\mathcal{F}}_{row}$ and build the corresponding filter bank $\{\phi_1^{row},\{\psi_n^{row}\}_{n=1}^{N_R}\}$,
\item Perform the boundaries detection on $\widetilde{\mathcal{F}}_{columns}$ and build the corresponding filter bank $\{\phi_1^{columns},\{\psi_n^{columns}\}_{n=1}^{N_C}\}$,
\item Filter $f$ along the rows with $\{\phi_1^{row},\{\psi_n^{row}\}_{n=1}^{N_R}\}$ which provides $N_R+1$ output images,
\item Filter each previous output image along the columns with $\{\phi_1^{columns},\{\psi_n^{columns}\}_{n=1}^{N_C}\}$, this provides at the end $(N_R+1)(N_C+1)$ subband images.\\
\end{enumerate}

\begin{figure}[!t]
\includegraphics[width=\columnwidth]{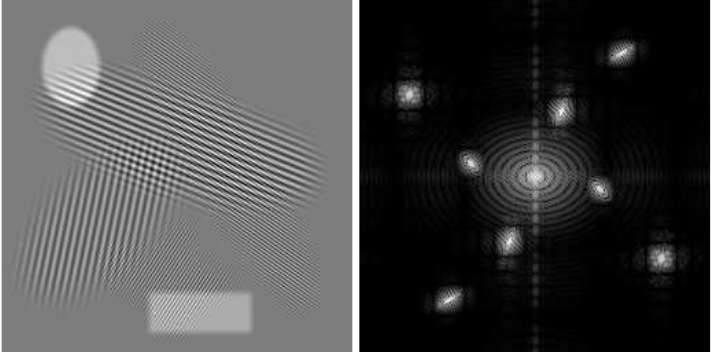}
\caption{Test image on left and the logarithm of its Fourier transform on right.}
\label{fig:im}
\end{figure}

\begin{figure}[!t]
\includegraphics[width=\columnwidth]{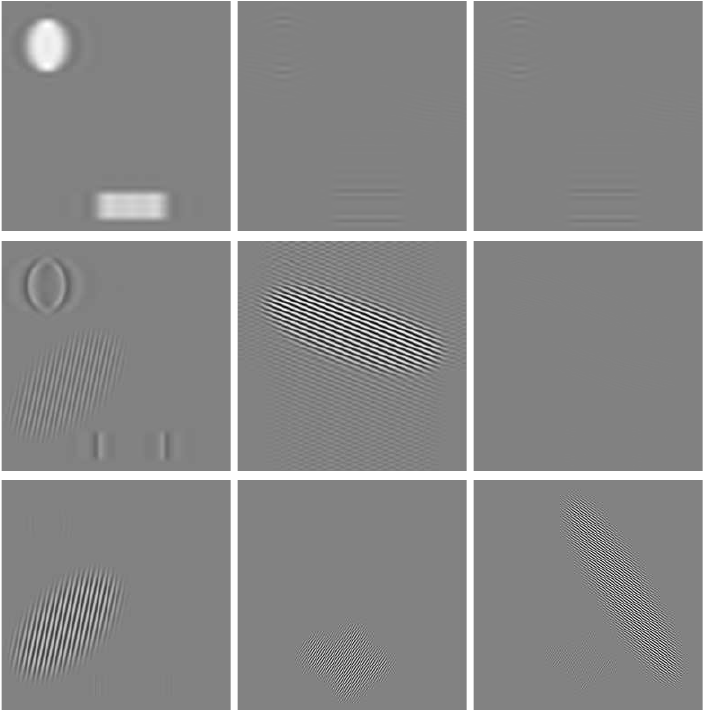}
\caption{Output of the 2D EWT of the test image.}
\label{fig:ewt2d}
\end{figure}

\begin{figure}[!t]
\centering\includegraphics[width=0.5\columnwidth]{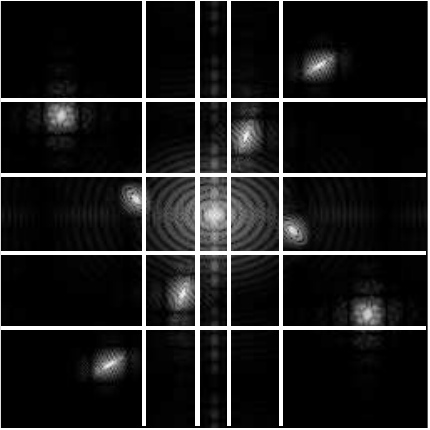}
\caption{Resulting partitioning of the 2D Fourier domain.}
\label{fig:2dtilling}
\end{figure}

In order to illustrate this 2D version of the EWT, we perform the previously described algorithm on the synthetic image given in left of Fig.~\ref{fig:im}. 
In order to see the adaptability of the 2D-EWT, this image is composed of objects and different 2D AM-FM components as shown in the Fourier domain 
(right of Fig.~\ref{fig:im}). We fix the maximum number of band per direction to two. The resulting subband images are given in Fig.~\ref{fig:ewt2d} and 
the corresponding empirical tilling of the Fourier domain in Fig.~\ref{fig:2dtilling}. We can see that the algorithm is capable to isolate each AM-FM 
component in different Fourier domains as we expected.

\section{Conclusion}\label{sec:conc}
In this paper, we proposed a novel approach to build wavelets adapted to represent the processed signal. The key idea is to build a wavelet filter bank based on Fourier supports detected from 
the information contained in the processed signal spectrum. We showed that it is possible to build a tight frame set of wavelets and then defined the Empirical Wavelet Transform (EWT) and its inverse. 
Like the classic wavelets, the proposed empirical wavelets correspond, in the temporal domain, to dilated versions of a single mother wavelet. The new insight lies in the fact the corresponding dilation factors don't 
follow a prescribed scheme but are detected empirically.
Some experiments made 
on simulated and real signals showed that this approach is able to separate different AM-FM modes composing the input signal. A comparison with the Empirical Mode Decomposition (EMD) showed that the EWT gives a more consistent 
decomposition while, generally, the EMD exhibits too much modes, which are sometime really difficult to interpret. Another advantage of the EWT compared to the EMD is that we can adapt the classic wavelet formalism to 
understand it. We also proposed a 2D extension based on the tensor product idea to process images. In this case, the 2D-EWT provides an adaptive tilling of the 
2D Fourier plane.\\
However, in some cases, the proposed method can fail (like any other wavelet approaches) compared to the EMD. For instance, if the input signal is composed of two chirps which overlap in both the time 
and frequency domains then the EWT will not be able to separate them while the EMD is supposed to be able to extract the most oscillating part first and the lowest one next. Such cases should probably be addressed 
by building adaptive frames with enough redundancy.\\
We also want to mention the differences with local cosine (sine) bases \cite{Auscher1992,Coifman1992}. Local cosine bases use a similar idea in their construction by splitting the time axis into 
several subintervals. Then a set of cosine wave forms, modulated by a bell function defined on each subinterval, is built and provides a basis of $L^2$ over the considered 
segment. If the general theory considers arbitrary segments, in practice it is too difficult to empirically determine the best segmentation and the existing 
algorithms return to the use of a tree structure based on prescribed dyadic subdivisions. Differences with the proposed approach are in the sense that we perform a real 
empirical segmentation in the Fourier domain instead of the temporal domain and we don't use a set of predefined functions.\\
In the future, we plan to investigate many directions based on the EWT. First, a full study about how to perform spectrum segmentation is necessary to design the best way to extract the different modes. Indeed as shown in 
the decomposition of $f_{Sig3}$ in the experiment section, taking the center between two local maxima does not use the information of the spectrum shape and return ``perturbed'' modes. 
Concerning the automatic detection of the number of modes in the signal, it could be interesting to use the concept of best basis like in the wavelet packets transform. This question can be also related to the previous 
one (how to detect the modes) if we take a clustering point of view where the estimation of $N$ corresponds to the estimation of the number of clusters. It will also be of interest to consider 
other kind of filters based on independent supports (like splines) and study their properties. \\
If the tensor product idea to generalize the proposed method to signals of any dimension arise naturally, we are extending the empirical adaptivity to 
other approaches like the ridgelet and curvelet transforms. The most general solution should be to directly perform a segmentation of the Fourier spectrum but 
this leads us to difficult and still open questions we expect to address in the future.\\
Another important direction to explore is the usefulness of this adaptive representation for applications like denoising, deconvolution, \ldots. These questions will be addressed in a forthcoming paper.\\
We also want to mention that a Matlab toolbox to perform the EWT will be available soon on the Matlab Central website.

\section{Acknowledgements}
The author wants to thank Prof. Andrea Bertozzi and Prof. Stanley Osher for their support and the fruitful discussions about this work. The author 
also wants to thanks the reviewers for their valuable comments which permit to improve the presentation of this work.



\begin{IEEEbiography}[{\includegraphics[width=1in,height=1.25in,clip,keepaspectratio]{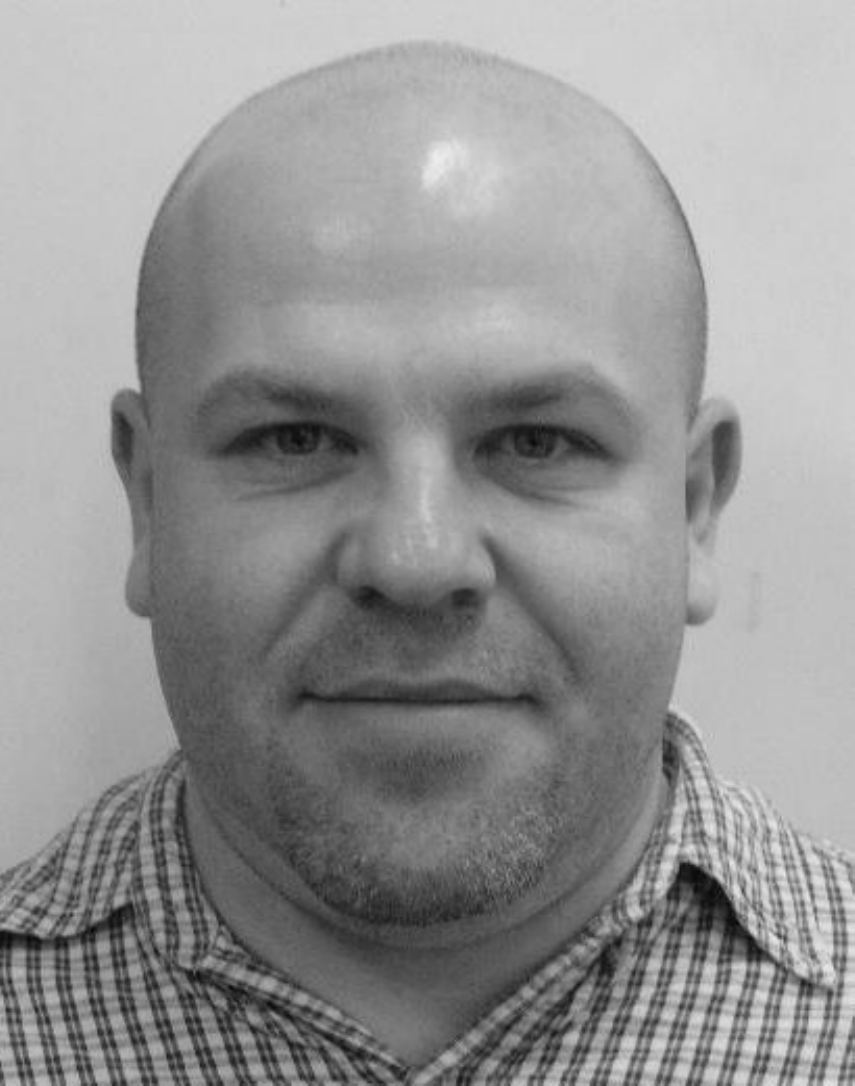}}]{J\'er\^ome Gilles}
J\'er\^ome Gilles is Assistant Adjunct Professor in the Department of Mathematics at UCLA in the team of Professor Stanley Osher.  
He received his Ph.D degree in mathematics (advised by Prof. Yves Meyer, Gauss Prize in 2010) from the Ecole Normale Sup\'erieure (ENS) of 
Cachan in France in 2006. Prior its position at UCLA, from 2001 to 2010, he held a position in the French Ministry of Defense in the Department of 
Space, Observation, Intelligence and UAVs systems, as an expert in signal and image processing. He also is a permanent associate research member 
of the Center of Mathematics and Their Applications (CMLA) of ENS Cachan in France.\\
His research interests are applied mathematics for signal/image modeling and analysis, textures analysis, signal/image restoration and the analysis 
and correction of atmospheric turbulence effects in long range imaging. He also devotes an important part of his time to study mathematical theories 
such as functional analysis, wavelets, partial differential equations and their applications for signal and image processing.\\
He wrote several international papers for SPIE conferences, Journal of Mathematical Imaging and Vision (JMIV), IEEE Trans. On Image Processing, 
Inverse Problems and Imaging Journal and had participated in many workshops and seminars in mathematics and computer vision. He also is a reviewer 
for many journals and international conferences.
\end{IEEEbiography}

\end{document}